\documentclass[12pt,reqno]{article}

\usepackage[usenames]{color}
\usepackage{amssymb}
\usepackage{graphicx}
\usepackage{amscd}

\usepackage[colorlinks=true,
linkcolor=webgreen,
filecolor=webbrown,
citecolor=webgreen]{hyperref}

\definecolor{webgreen}{rgb}{0,.5,0}
\definecolor{webbrown}{rgb}{.6,0,0}

\usepackage{color}
\usepackage{fullpage}
\usepackage{float}

\usepackage{graphics,amsmath,amssymb}
\usepackage{amsthm}
\usepackage{amsfonts}
\usepackage{latexsym}
\usepackage{epsf}

\makeatletter

\theoremstyle{plain}
\newtheorem{theorem}{Theorem}
\newtheorem{lemma}[theorem]{Lemma}

\theoremstyle{definition}

\newcommand{\summ}{\sum\limits}

\theoremstyle{definition}

\makeatother

\begin{document}

\begin{center}
\vskip 1cm{\LARGE\bf
The central component of a triangulation
}\\
\large
Alon Regev\\
Department of Mathematical Sciences\\
Northern Illinois Univeristy\\
DeKalb, IL\\
\end{center}

\section{Introduction}
\begin{figure}
\begin{center}
\epsfxsize=1.6in
\leavevmode\epsffile{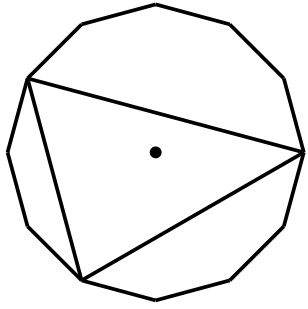}
\end{center}
\caption{A central triangle with $n=12$, $i=3$, $j=4$ and $k=5$.}\label{Fig1}
\end{figure}

	Considering a triangulation of a regular convex polygon as a subset of $\mathbb{R}^2$ centered at the origin, define its {\em central component} to be the diameter or triangle that contains the origin (see Figure \ref{Fig1}). More generally, every dissection of a polygon can be associated with its set of components, including one central component. To our knowledge, the notion of components and central components was first used by  D. Bowman and the author to study symmetry classes of dissections \cite{BR2}. In this note we use central components to obtain new recursion relations for Catalan and $k$-Catalan (Fuss-Catalan) numbers, and use these recursions to prove congruence relations of these numbers. We also enumerate the triangulations that include a fixed vertex in their central components.


	Let $C_n$ be the $n$-th Catalan number, so $C_{n-2}$ is the number of triangulations of an $n$-gon, and let $C_x=0$ unless $x$ is an integer. The Catalan numbers satisfy the recursion relation
\begin{align}\label{standard}
\summ_{k=0}^nC_k	C_{n-k}=C_{n+1}.
\end{align}

Consider a triangulation of an $n$-gon as a labeled graph with vertices $0,1,\ldots, n-1$ and edges denoted $xy$ for distinct vertices $x$ and $y$. The edges include $n$ sides $01, 12, \ldots ,(n-1)0$ and $n-3$ diagonals. The {\em cyclic length} of an edge $xy$, with $x<y$, is defined as \[\min \{y-x, n+x-y\}.\]

By enumerating the triangulations of an $n$-gon according to their central components, we obtain the following recursion relation.
\begin{lemma}\label{center-lem}
\begin{align}\label{center-eq}
C_{n-2}={n\over 2} C_{n/2-1}^2+\summ_{\substack{i+j+k=n\\ i\le j \le k<n/2}}m_{ijk}C_{i-1}C_{j-1}C_{k-1},
\end{align}
where
\begin{align*}
m_{ijk}=
\begin{cases}
{n\over 3} & \text{if }i=j=k;\\
n  & \text{if }i<j=k \text{ or } i=j<k;\\
2n & \text{if } i<j<k.\\
\end{cases}
\end{align*}
\begin{proof}
The first term enumerates the triangulations whose central component is a diameter. In the summation, $m_{ijk}$ is the number of ways to place a triangle whose sides have cyclic lengths lengths $i, j, k$ inside an $n$-gon (see Figure \ref{Fig1}). The conditions under the summation ensure that indeed this is a central triangle.
The three cases determining $m_{ijk}$ correspond to the central triangle being equilateral, isosceles or scalene.
\end{proof}
\end{lemma}


\section{Congruences relations}

Congruence relations of Catalan numbers $C_n$ and related sequences have been the object of extensive study (see \cite{AK,DS} and the references therein). We next show that Lemma \ref{center-lem} can be used to derive some results of this nature.

Note that the following result can be proved more easily using \eqref{standard}.
\begin{theorem}\label{oddcats}
$C_n$ is odd if and only if $n=2^a-1$ for some integer $a\ge 0$.
\begin{proof}
We use induction on $n$, with the base cases easily verified.

Let $n=2^a+1$, and we use \eqref{center-eq} show that $C_{n-2}$ is odd. Now $C_{n/2-1}=0$, and we show that there is exactly one odd term in the summation, namely
\begin{align}\label{ot}
m_{1, 2^{a-1}, 2^{a-1}}C_0C_{2^{a-1}-1}^2=n C_{2^{a-1}-1}^2.
\end{align}

By induction this term is odd, so it remains to show there are no other odd terms.
Again by induction, any odd term has $i=2^b, j=2^c, k=2^d$ for some $0\le b\le c\le d$.
Therefore
$2^b+2^c+2^d=2^a+1$, which is only possible when $b=0$ and $c=d=a-1$. Thus we have proved that $C_n$ is odd when $n=2^a-1$.

Now suppose $C_{n-2}$ is odd.  If $n$ is even, then all of the terms in the summation are even, so that ${n\over 2} C_{n/2-1}^2$ must be odd.
Therefore $n/2$ is odd, and at the same time by induction $n/2-1=2^b-1$ for some $b\ge 0$. Thus $b=0$ and $n=0=2^0-1$.

If $n$ is odd then $C_{n/2-1}=0$. Any odd term in the summation has one of $i, j, k$ odd, but by induction $i-1,j-1,k-1$ all have the form $2^b-1$.
It follows that $i=1$, and the odd term is the one given by \eqref{ot}.

\end{proof}
\end{theorem}

We next use Lemma \ref{center-lem} to prove a recent result of Eu, Liu and Yeh \cite{ELY}.
\begin{theorem}\cite[Theorem 2.3]{ELY}
For all $n\ge 0$,
\begin{align}\label{mod4eq}
C_n\equiv_4
\begin{cases}
1 & \text{ if } n=2^a-1 \text{ for some } a\ge 0;\\
2 & \text{ if } n=2^a+2^b-1 \text{ for some } b>a\ge 0;\\
0 & \text{ otherwise.}\\
\end{cases}
\end{align}
\begin{proof}
By verifying the base cases we may assume $n\ge 10$ and proceed by induction on $n$. Consider the three cases on the right hand side of \eqref{mod4eq}.

\begin{enumerate}
\item The case $n=2^a-1$:

 Let $n=2^a+1$, and we show that $C_{n-2}\equiv_4 1$. Reducing \eqref{center-eq} modulo $4$ gives
\begin{equation}\label{mod4}
C_{n-2}\equiv_4 \summ_{\substack{i+j+k=n\\ i\le j \le k<n/2}}m_{ijk}C_{i-1}C_{j-1}C_{k-1}.
\end{equation}
By Theorem \ref{oddcats},
\begin{equation}\label{t1}
nC_{2^{a-1}}C_{2^{a-1}}C_0\equiv_4 1.
\end{equation}
It remains to show that this is the only term not divisible by $4$.

If $i=j=k$ and $C_{i-1}$ is odd then $i=2^c$ for some $c\ge 0$, so $3(2^c-1)=2^a+1$. Thus $2^{c+1}+2^c=2^a+2^2$, which is impossible for $n\ge 10$.

Suppose $i<j<k$ with $C_{i-1}C_{j-1}C_{k-1}$ odd. Then $i=2^c$, $j=2^d$ and $k=2^e$ for some $e>d>c\ge 0$, but then $2^e+2^d+2^c=2^a+1$, which is impossible. Thus
$2nC_{i-1}C_{j-1}C_{k-1}\equiv_4 0$ in this case.

Now suppose $i=j<k$. If $C_{i-1}C_{j-1}C_{k-1}\not \equiv_4 0$ then by induction, $i=2^c$ for some $c\ge 0$ and either $k=2^d$ with $d>c$, or $k=2^d+2^e$ with $e>d\ge c$.
Assume $k=2^d$. Then $2^a+1=2\cdot 2^c+2^d$, so $c=0$ and $d=a$, yielding the term given by \eqref{t1}. If $k=2^d+2^e$ then $2^a+1=2\cdot 2^c+2^d+2^e$, a contradiction.
Similarly, assuming $i<j=k$ and $C_{i-1}C_{j-1}C_{k-1}\not\equiv_4 0$ results in a contradiction.

\item The case $n=2^a+2^b-1$:

Let $n=2^a+2^b+1$, and we show that $C_{n-2}\equiv_4 2$. By verifying the base cases we may assume $b>a\ge 1$.
Note that \eqref{mod4} holds in this case as well. Now
\begin{equation}\label{2t}
nC_{2^a-1}C_{2^{b-1}-1}C_{2^{b-1}-1}\equiv_4 2.
\end{equation}
 We show that all other terms are divisible by $4$.
Since the binary expansion of $n$ is unique, by Theorem \ref{oddcats} there are no terms with $C_{i-1}C_{j-1}C_{k-1}$ odd.
Therefore the terms $2nC_{i-1}C_{j-1}C_{k-1}$ with $i<j<k$ are divisible by $4$ and the term with $i=j=k$ is divisible by $8$.
It remains to show that the term given in \eqref{2t} is the only one with $i=j<k$ or $i<j=k$ that is congruent to $2$.

Suppose then that $2^c+2^d+2\cdot 2^e = 2^a+2^b+1$ for some $d>c\ge 0$ and $e\ge 0$. Then $c=0$ and $2^d+2^{e+1}=2^a+2^b$.
Now $d=b$ would imply $2^c+2^d>n/2$. Therefore $d=a$ and $e=b-1$, as asserted.

\item {Otherwise:}

 We show that $C_{n-2}\equiv_40$ unless $n-2$  has one of the forms above.
If $C_{n/2-1}^2\not \equiv_4 0$ then by Theorem \ref{oddcats} we have $n/2-1=2^a-1$ for some $a\ge 0$. Therefore $n=2^{a+1}$. For $n$ sufficiently large this implies $8|n$  so that
${n\over 2}C_{n/2-1}^2\equiv_4 0$.

Next consider the terms in the summation.

\begin{enumerate}
\item {The terms $i=j=k$:}

 If $C_{n/3-1}^3\not \equiv_4 0$ then by Theorem \ref{oddcats} we have
$n/3-1=2^a-1$ for some $a\ge 0$, so $n=3\cdot 2^a$, which for $n$ sufficiently large implies that $n/3$ is divisible by $4$. Thus $m_{ijk}C_{i-1}C_{j-1}C_{k-1}\equiv_4 0$.

\item {The terms $i=j<k$ or $i<j=k$:}

Suppose $nC_{i-1}C_{j-1}C_{k-1}\not \equiv_4 0$.
If $C_{i-1}C_{j-1}C_{k-1}$ is odd then by Theorem \ref{oddcats} we have $n=2^{b+1}+2^c$ for some $b,c\ge 0$. Now $2^c<n/2=2^b+2^{c-1}$, so $c-1< b$ and in fact $c<b$.
Since $4\not | n$, this implies $c=0$ or $c=1$, so $n-2=2^{n+1}-1$ or $n-2=2^{b+1}+2^0-1$.

The only case left to check is when $n$ is odd and $C_{i-1}C_{j-1}C_{k-1}\equiv_4 2$. In this case, by induction $n=2\cdot 2^b+(2^c+2^d)$ with $b,c,d\ge 0$ and $d>c$.
Since $n$ is odd then $c=0$, so that $n-2=2^{b+1}+2^d-1$ as asserted.
\item {The terms $i<j<k$:}

If $2nC_iC_jC_k\ne 4$ then $n$, $C_i$, $C_j$ and $C_k$ are all odd. By Theorem \ref{oddcats} we have $n=2^b+2^c+2^d$ for $d>c>b\ge 0$, and since $n$ is odd then $b=1$ and $n-2=2^c+2^d-1$.
\end{enumerate}
\end{enumerate}
 \end{proof}
 \end{theorem}



Another congruence relation follows immediately from reducing \eqref{center-eq} modulo $p$.

\begin{theorem}\label{modpthm}
If $p\ge 5$ is prime and $n\equiv_p -2$ then $C_n\equiv_p 0$.
\end{theorem}

\section{Generalization to $k$-angulations}

Lemma \ref{center-lem} can be generalized to give a recursion for the number of {\em $k$-angulations}, which are partitions of a polygon into $k$-gons.
For example, let
\[Q_n={1 \over 2n+1}{3n\choose n}\]
be the number of quadrangulations of a $(2n+2)$-gon, and let $Q_x=0$ unless $x$ is a nonnegative integer.
\begin{lemma}
\begin{align}\label{quads-eq}
Q_n=(n+1)Q_{n/2}^2+\summ_{\substack{i+j+k+l=2n+2\\ i\le j \le k\le l<n+1}}m_{ijkl}Q_{(i-1)/2}Q_{(j-1)/2}Q_{(k-1)/2}Q_{(l-1)/2},
\end{align}
where
\begin{align}
m_{ijkl}=
\begin{cases}
{N\over 4}& \text{ if } i=l,\\
N& \text{ if } i=k<l \text{ or } i<j=l,\\
{3N\over 2}& \text{ if } i=j<k=l,\\
3N & \text{ if } i=j<k<l \text{ or } i<j=k<l \text{ or } i<j<k=l, \\
6N & \text{ if } i<j<k<l\\
\end{cases}
\end{align}
for $N=i+j+k+l$.
\end{lemma}

More generally, let
$f_{n,k}$ be the number of $k$-angulations of an $n$-gon. It is well known that

\begin{align}\label{kang-Fuss-eq}
 f_{(k-1)n+2,k+1}=C_{n,k}
 \end{align}

where
 \[C_{n,k}={1\over (k-1)n+1}{kn\choose n}\]
are the Fuss-Catalan (or $k$-Catalan) numbers.
Define $f_{n,k}=0$ unless $n=(k-2)m+2$ for some integer $m\ge 0$.

\begin{lemma}
For any $n\ge 2$ and $ k\ge 3$,
\begin{align}\label{kang-eq}
f_{n,k}={n\over 2}f_{n/2+1,k}+\summ_{\substack{i_1+\ldots+ i_k=n\\ i_1\le \ldots \le i_k <n/2}} m_{i_1\ldots i_k} f_{i_1+1,k} \cdots f_{i_k+1, k},
\end{align}
where $m_{i_1\ldots i_k}$ is the number of ways position a $k$-gon with sides of cyclic lengths $i_1,\ldots, i_k$ inside an $N$-gon for $N=i_1+\ldots+i_k$.
\end{lemma}

Theorem \ref{modpthm} can be generalized by reducing \eqref{kang-eq} modulo $p$.

\begin{theorem}
If $p\ge 3$ is prime with $p\not | k$ and $p | n$
then $f_{n,k}\equiv_p 0$.
\begin{proof}
For a given $k$-gon, the number of cyclic permutations of the $k$ sides that leave the $k$-gon unchanged is divisible by $k$. Therefore the number of inequivalent rotations of the
$k$-gon inside the $n$-gon is divisible by $n/k$. It follows that $m_{i_1\ldots i_k}$ is divisible by $n/ k$, and so the given assumptions imply $p$ that divides $f_{n,k}$.
\end{proof}
\end{theorem}


\section{Triangulations with a fixed vertex in their central component}

L. Shapiro 
proposed the following question: how many triangulations include the vertex $0$ in their central component? The following theorem answers this question.

\begin{figure}
\begin{center}
\epsfxsize=1.8in
\leavevmode\epsffile{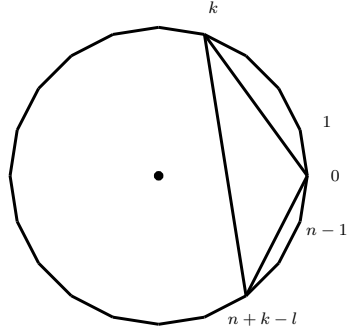}
\end{center}
\caption{Triangulations with $0$ outside the central component.}\label{Fig2}
\end{figure}

\begin{theorem}
The number $a(n)$ of triangulations of an $n$-gon with the vertex $0$ outside their central component is
\[a(n)=\summ_{m=1}^{\lfloor n/2 \rfloor-1}C_mC_{n-2-m}.\]
\begin{proof}
Enumerate these triangulations according to the cyclic length $l$ of the shortest diagonal that separates $0$ form the center (see Figure \ref{Fig2}).
Note that $2\le l\le \lfloor n/2 \rfloor$.  Given such $l$, suppose this diagonal is given by $k(n+k-l)$. Note that $1\le j\le l-1$.
Since this is the shortest such diagonal, the triangulation must also
include the diagonals $0k$ and $0(n+k-l)$, forming a triangle. The regions outside of this triangle can be triangulated arbitrarily. Therefore
\begin{align*}
a(n)&=\summ_{l=2}^{\lfloor n/2 \rfloor} \summ_{k=1}^{l-1} C_{n-l-1}C_{l-k-1}C_{k-1}\\
	&=\summ_{l=2}^{\lfloor n/2 \rfloor}  C_{n-l-1}C_l\\
	&=\summ_{m=1}^{\lfloor n/2 \rfloor-1}C_mC_{n-2-m},\\
\end{align*}
where the second equality follows from \eqref{standard}.
\end{proof}
\end{theorem}

It seems that the sequence $a(n)$ is given by  \cite[A027302]{Sl},
which is the number of Dyck paths with $UU$ spanning their midpoint, and a formula is given there by:

\[a(n)=\summ_{0\le k < n/2}T(n,k)T(n,k+1),\]

where \[T(n,k)=  {n-2k+1\over n-k+1}{n \choose k}.\]

It would be interesting to find a bijection between these sets of triangulations and Dyck paths.

\end{document}